\documentclass[12pt,a4paper]{article}
\usepackage[english]{babel}
\usepackage[utf8]{inputenc}

\usepackage{amsmath}
\usepackage{amssymb}
\usepackage{amsfonts}
\usepackage{amsbsy}

\usepackage{graphicx}

\usepackage{newlfont}
\usepackage{float}
\usepackage[paperwidth=192mm,
            		   paperheight=262mm,
           		 	   vmargin={19mm,19mm},
            		   hmargin={13.7mm,13.7mm},
            		   headsep=12pt,
            		   footskip=12pt]{geometry}

\usepackage{hyperref}

\usepackage{siunitx}
\usepackage{subcaption}

\usepackage{caption}
\usepackage{multirow}

\usepackage{booktabs}
\usepackage{tabularx}

\usepackage{placeins}

\usepackage[normalem]{ulem}

\usepackage[backend=bibtex,style=numeric,sorting=nyt]{biblatex}
\renewbibmacro{in:}{}
\DeclareMathOperator{\grad}{\nabla}
\DeclareMathOperator{\dive}{\nabla\cdot}

\begin{document}

\title{Efficient and scalable atmospheric dynamics simulations using non-conforming meshes}

\author{Giuseppe Orlando$^{(1)}$, \\ 
        	 Tommaso Benacchio$^{(2)}$, Luca Bonaventura$^{(3)}$}

\date{}

\maketitle

\begin{center}
{
\small
$^{(1)}$  
CMAP, CNRS, \'{E}cole polytechnique, Institut Polytechnique de Paris \\ Route de Saclay, 91120 Palaiseau, France \\
{\tt giuseppe.orlando@polytechnique.edu} \\
\ \\
$^{(2)}$  
Weather Research, Danish Meteorological Institute \\
Sankt Kjelds Plads 11, 2100 Copenhagen, Denmark \\
{\tt tbo@dmi.dk} \\
\ \\
$^{(3)}$  
Dipartimento di Matematica, Politecnico di Milano \\
Piazza Leonardo da Vinci 32, 20133 Milano, Italy \\
{\tt luca.bonaventura@polimi.it} \\
}
\end{center}

\noindent


\pagebreak

\abstract{We present the massively parallel performance of a $h$-adaptive solver for atmosphere dynamics that allows for non-conforming mesh refinement. The numerical method is based on a Discontinuous Galerkin (DG) spatial discretization, highly scalable thanks to its data locality properties, and on a second order Implicit-Explicit Runge-Kutta (IMEX-RK) method for time discretization, particularly well suited for low Mach number flows. Simulations with non-conforming meshes for flows over orography can increase the accuracy of the local flow description without affecting the larger scales, which can be solved on coarser meshes. We show that the local refining procedure has no significant impact on the parallel performance and, therefore, both efficiency and scalability can be achieved in this framework.}

\pagebreak

\section{Introduction}
\label{sec:intro}

Efficient numerical simulations of atmospheric flows are crucial for weather and climate predictions and pose several computational challenges. Peculiar to this application are the timeliness constraints that dictate maximum admissible simulation runtimes in operational weather prediction. Medium range numerical weather prediction (NWP) forecasts up to ten days ahead are typically expected to complete within one hour in the forecast cycles of weather centres, thus imposing demanding modelling choices in order to guarantee computational efficiency. In a context of increasing need for computational resources due to higher spatial resolutions, massively parallel model scalability is therefore required.

From a physical standpoint, slow-moving atmospheric flows of meteorological interest are characterized by a speed much lower than the speed of sound, so that their Mach number is low and compressibility effects are usually deemed not very relevant. Weakly compressible flows are an example of problem with multiple length and time scales \cite{klein:1995}. Moreover, atmospheric flows display phenomena on a very wide range of spatial scales that interact with each other. Strongly localized features require a very high spatial resolution to be correctly resolved, while larger scale features, such as midlatitude pressure systems and stratospheric flows, can be adequately resolved on coarser meshes. Hence, the design of efficient and stable numerical schemes for such models is a challenging task. 

Because of its multi-scale nature, NWP and, in particular, flows over orography are an apparently ideal framework to develop adaptive numerical approaches based on variable resolution meshes. However, mesh adaptation strategies have slowly found their way into the NWP literature, due to concerns about the accuracy of variable resolution meshes for the correct representation of atmospheric wave phenomena, and the greater complexity of an efficient parallel implementation for non-uniform or adaptive meshes. Moreover, numerical strategies with variable resolution meshes typically employ local mesh refinement only in the horizontal directions, while columns of cells with the same horizontal dimension are employed in the vertical direction \cite{kopera:2014, muller:2013, yelash:2014}. A non-conforming mesh is characterized by neighbouring cells with different resolution on both the horizontal and the vertical direction \cite{orlando:2024a}. In \cite{hellsten:2021}, a full 3D nesting approach for atmospheric flows is presented. However, the method is tested only on cases without orography. To the best of our knowledge, in \cite{orlando:2024a} the authors firstly proposed a method able to decrease both horizontal and vertical resolution as height increases, filling a gap in the NWP literature and showing how fully 3D non-conforming meshes can be successfully employed for flows over orography. The solver is based on the IMEX-DG method proposed in \cite{orlando:2022} (see also \cite{orlando:2024b}) and employed for atmospheric flows in \cite{orlando:2023, orlando:2024c, orlando:2024a}. Thanks to its data locality features, DG simulations are characterized by small communication-to-computation ratios and increasingly good scalability at higher orders of accuracy.   

In this work, we present the parallel performance of the solver implementation, carried out in the framework of the \texttt{deal.II} library \cite{arndt:2023, bangerth:2007} which natively allows for the use of non-conforming meshes. Using a well-established library with an active user community allows the investigation of advanced numerical choices without the need to code basic features of the numerical method or to implement parallel paradigms. Here we show that the local mesh refinement procedure does not adversely affect the parallel efficiency and scalability of the model compared to the simulations using uniform meshes, as measured in runs performed on the MeluXina EuroHPC high-performance computing facility.

The manuscript is structured as follows. In Section \ref{sec:model_num}, we briefly review the model equations and the numerical framework. The application and the analysis of parallel performance for a relevant benchmark is presented in Section \ref{sec:test}. Finally, some conclusions are reported in Section \ref{sec:conclu}. 

\section{The model equations and the numerical framework}
\label{sec:model_num}

The compressible Euler equations of gas dynamics represent the most comprehensive mathematical model for atmosphere dynamics \cite{giraldo:2008, steppeler:2003}. Let $\Omega \subset \mathbb{R}^{d}, 2 \le d \le 3$ be a domain and denote by $\mathbf{x}$ and $t$ the spatial coordinates and the temporal coordinate, respectively. The mathematical model reads as follows:
\begin{eqnarray}\label{eq:euler_comp}
	\frac{\partial\rho}{\partial t} + \dive\left(\rho\mathbf{u}\right) &=& 0 \nonumber \\
	\frac{\partial\left(\rho\mathbf{u}\right)}{\partial t} + \dive\left(\rho\mathbf{u} \otimes \mathbf{u}\right) + \grad p &=& \rho\mathbf{g} \\
	\frac{\partial\left(\rho E\right)}{\partial t} + \dive\left[\left(\rho E + p\right)\mathbf{u}\right] &=& \rho \mathbf{g} \cdot \mathbf{u}, \nonumber
\end{eqnarray}
for $\mathbf{x} \in \Omega, t \in (0, T_{f}]$, supplied with suitable initial and boundary conditions. Here $T_{f}$ is the final time, $\rho$ is the density, $\mathbf{u}$ is the fluid velocity, $p$ is the pressure, and $\otimes$ denotes the tensor product. Moreover, $\mathbf{g} = -g\mathbf{k}$ is the acceleration of gravity, with $g = \SI{9.81}{\meter\per\second\squared}$ and $\mathbf{k}$ being the upward pointing unit vector in the standard Cartesian frame of reference. Finally, $E$ denotes the total energy per unit of mass and we point out the one can rewrite $\rho E = \rho e + \rho k$, where $e$ is the internal energy and $k = \frac{1}{2}\|\mathbf{u}\|^{2}$ is the kinetic energy. System \eqref{eq:euler_comp} is complemented by the equation of state of ideal gases, given by $p = \rho RT$, where $R$ is the specific gas constant and $T$ denotes the temperature. We take $R = \SI{287}{\joule\per\kilo\gram\per\kelvin}$. 

\noindent

A dimensional analysis can be carried out, leading to the following system of equations \cite{klein:1995, munz:2003, orlando:2022, orlando:2024b}:
\begin{eqnarray}\label{eq:euler_comp_adim}
	\frac{\partial\rho}{\partial t} + \dive\left(\rho\mathbf{u}\right) &=& 0 \nonumber \\
	\frac{\partial\left(\rho\mathbf{u}\right)}{\partial t} + \dive\left(\rho\mathbf{u} \otimes \mathbf{u}\right) + \frac{1}{M^{2}}\grad p &=& -\frac{1}{Fr^{2}}\rho\mathbf{k} \\
	\frac{\partial\left(\rho E\right)}{\partial t} + \dive\left[\left(\rho e + M^{2}\rho k + p\right)\mathbf{u}\right] &=& -\frac{M^{2}}{Fr^{2}}\rho\mathbf{k} \cdot \mathbf{u}, \nonumber
\end{eqnarray}
where, with a slight abuse of notation, the non-dimensional variables have the same symbols of the dimensional ones. Here, $M$ denotes the Mach number, i.e. the ratio between the local fluid speed and the speed of sound, while $Fr$ is the Froude number, i.e. the ratio between the flow inertia and gravitational forcing. 

Atmospheric flows, as those considered in this work, are characterized by low Mach number values. In the low Mach number limit, pressure gradient terms yield stiff components for the resulting semi-discretized ODE
system, since the pressure gradients in \eqref{eq:euler_comp_adim} are proportional to $\frac{1}{M^{2}}$. Hence, following \cite{casulli:1984, dumbser:2016}, the method proposed in \cite{orlando:2022} couples implicitly the energy equation to the momentum equation, while the continuity equation is treated in a fully explicit fashion. High-order accuracy in time is then achieved making use of Implicit-Explicit Runge-Kutta (IMEX-RK) time integrators \cite{kennedy:2003}, which are widely employed for ODE systems that include both stiff and non-stiff components. Finally, for the spatial discretization we employ the Discontinuous Galerkin (DG) method \cite{giraldo:2020}, which combines high-order accuracy and flexibility in a highly data-local framework. More specifically, we aim at employing non-conforming meshes, i.e. meshes for which the resolution between two neighbouring cells can be different along both horizontal and vertical direction. The DG method naturally allows the use of this kind of meshes \cite{heinz:2023} without any hanging node appearing. We refer to \cite{orlando:2024a} for a short introduction to non-conforming meshes, and to \cite{orlando:2022, orlando:2023, orlando:2024b} for a complete analysis and discussion of the numerical methodology.

\section{Numerical results}
\label{sec:test}

In this Section, we consider an idealized three-dimensional test case of an atmospheric flow over orography \cite{lock:2012, melvin:2019, orlando:2024a}. Simulation parameters are related to two Courant numbers, the so-called acoustic Courant number $C$, which is based on the speed of sound $c$, and the advective Courant number $C_{u}$, which is based on the speed of the local flow velocity $u$:
$$ C = rc\Delta t/\mathcal{H}\sqrt{d} \qquad C_{u} = ru\Delta t/\mathcal{H}\sqrt{d}. $$
Here, $r$ is the polynomial degree used for the DG spatial discretization, $\mathcal{H}$ is the minimum cell diameter of the computational mesh, and $\Delta t$ is the time step adopted for the time discretization. We consider polynomial degree $r = 4$. In the following, we analyze the accuracy of the simulations with mesh refinement and then focus on the scalability of the numerical model. The 9.5.2 \texttt{deal.II} \cite{arndt:2023, bangerth:2007} release has been used to produce the results in this section. The simulations have been run using up to 1024 2x AMD EPYC Rome 7H12 64c\@ 2.6GHz CPUs at MeluXina supercomputer \footnote{\url{https://docs.lxp.lu/}} and OpenMPI 4.1.5 has been employed.

\subsection{3D medium-steep bell-shaped hill}

We consider a three-dimensional configuration of a flow over a bell-shaped hill, originally proposed in \cite{lock:2012} and also employed in \cite{melvin:2019, orlando:2023, orlando:2024a}. The computational domain is $\Omega = \SI[parse-numbers=false]{\left(0, 60\right) \times \left(0, 40\right) \times \left(0, 16\right)}{\kilo\meter}$. The mountain profile is defined as follows:
\begin{equation}
	h(x,y) = \frac{h_{c}}{\left[1 + \left(\frac{x - x_{c}}{a_{c}}\right)^{2} + \left(\frac{y - y_{c}}{a_{c}}\right)^{2}\right]^{\frac{3}{2}}},
\end{equation}
with $h_{c} = \SI{400}{\meter}, a_{c} = \SI{1}{\kilo\meter}, x_{c} = \SI{30}{\kilo\meter},$ and $y_{c} = \SI{20}{\kilo\meter}$. The buoyancy frequency is $N = \SI{0.01}{\per\second}$, whereas the background velocity is $\overline{u} = \SI{10}{\meter\per\second}$. The final time is $T_{f} = \SI{1}{\hour}$. The initial conditions read as follows \cite{benacchio:2014}:
\begin{eqnarray}
	p &=& p_{ref}\left\{1 - \frac{g}{N^{2}}\Gamma\frac{\rho_{ref}}{p_{ref}}\left[1 - \exp\left(-\frac{N^{2}z}{g}\right)\right]\right\}^{1/\Gamma} \\
	\rho &=& \rho_{ref}\left(\frac{p}{p_{ref}}\right)^{1/\gamma}\exp\left(-\frac{N^{2}z}{g}\right),
\end{eqnarray}
where $p_{ref} = \SI[parse-numbers = false]{10^{5}}{\pascal}$ and $\rho_{ref} = \frac{p_{ref}}{R T_{ref}}$, with $T_{ref} = \SI{293.15}{\kelvin}$. Finally, we set $\Gamma = \frac{\gamma - 1}{\gamma}$, with $\gamma = 1.4$ being the isentropic exponent. Wall boundary conditions are employed for the bottom boundary, whereas non-reflecting
boundary conditions are required by the top boundary and the lateral boundaries. We refer to \cite{orlando:2024a} for the implementation of non-reflecting boundary conditions.

\begin{figure}[h!]
	\centering
	\includegraphics[width = 0.8\textwidth]{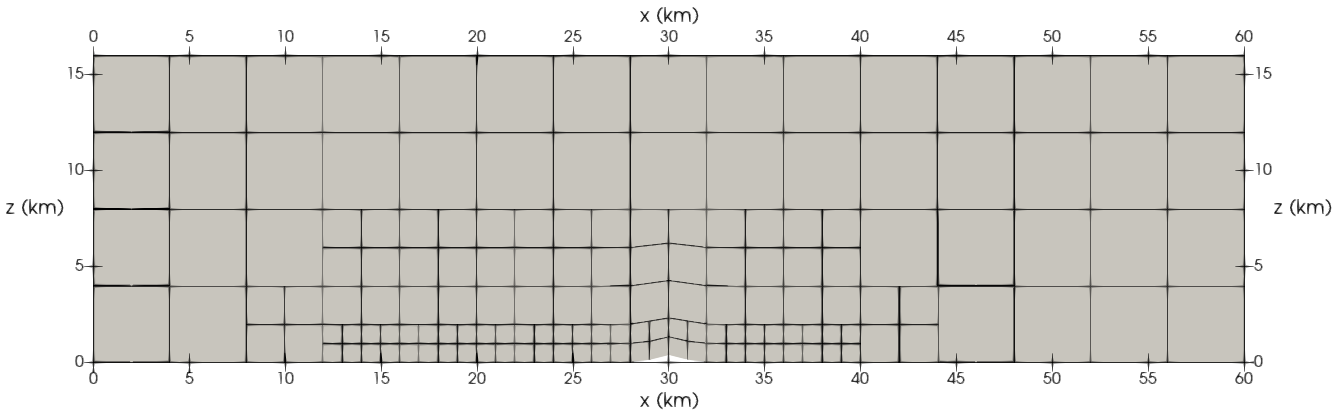}
	\caption{3D medium-steep bell-shaped hill test case, non conforming mesh. $x-z$ slice at $y = \SI{20}{\kilo\meter}$.}
	\label{fig:3D_non_conforming_mesh}
\end{figure}

We consider two different meshes: a uniform mesh composed by $60 \times 40 \times 16$ elements, i.e. a spatial resolution of $\SI{250}{\meter}$ along all the directions, and a non-conforming mesh composed by $N_{el} = 1958$, with its finest resolution corresponding to that of the uniform mesh (Figure \ref{fig:3D_non_conforming_mesh}). Notice that the resolution depends on the polynomial degree $r$ employed for the spatial discretization. More specifically, the effective resolution is computed dividing the size of the element along each direction by the polynomial degree. We take $\Delta t = \SI{2}{\second}$, yielding a maximum acoustic Courant number $C \approx 2.75$ and a maximum advective Courant number $C_{u} \approx 0.13$ for the finest uniform mesh. 

\begin{figure}[h!]
	\centering
	\begin{subfigure}{0.75\textwidth}
		\centering
		\includegraphics[width = 0.925\textwidth]  {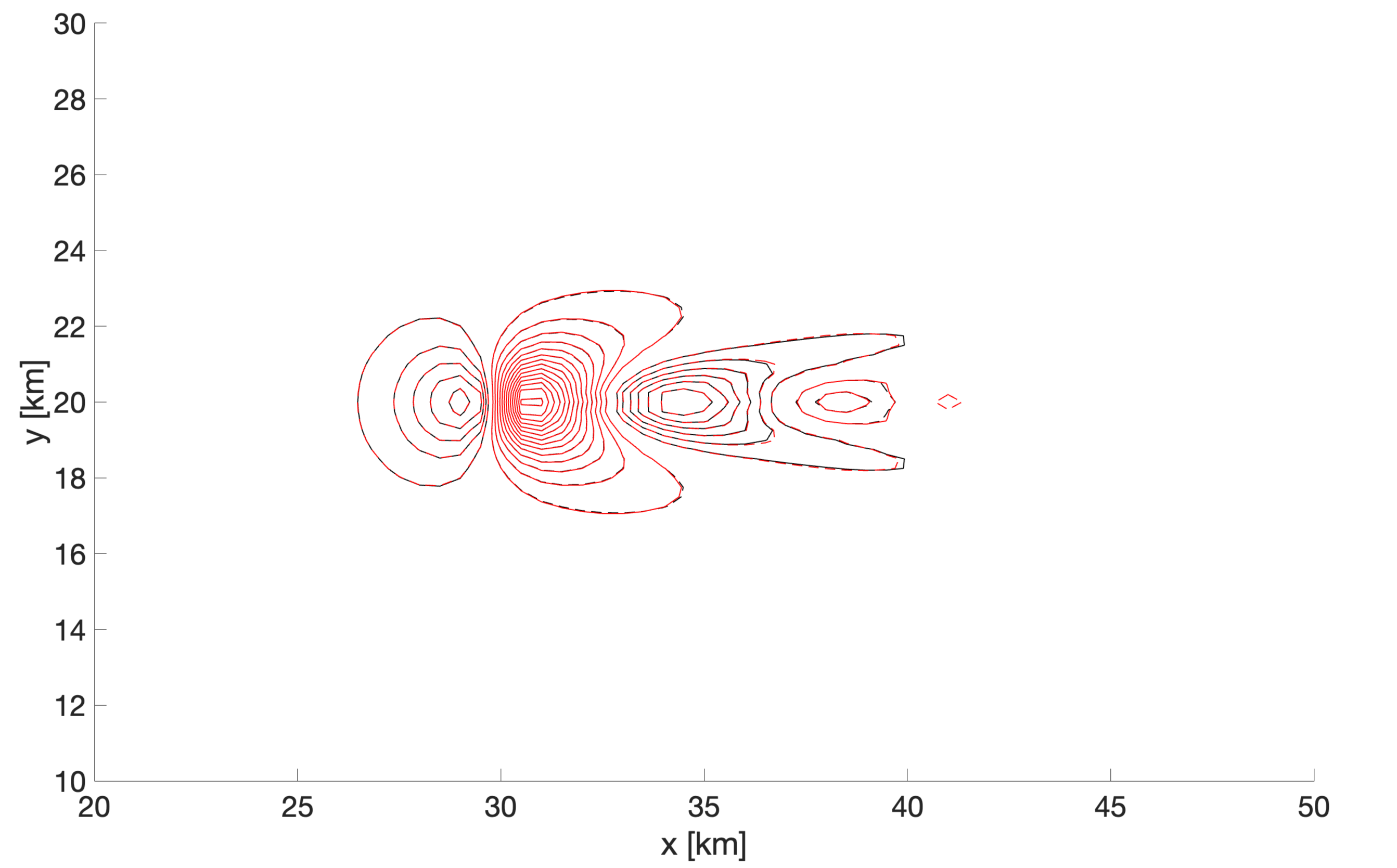}
	\end{subfigure}
	\\
	\begin{subfigure}{0.75\textwidth}
		\centering
		\includegraphics[width = 0.925\textwidth]{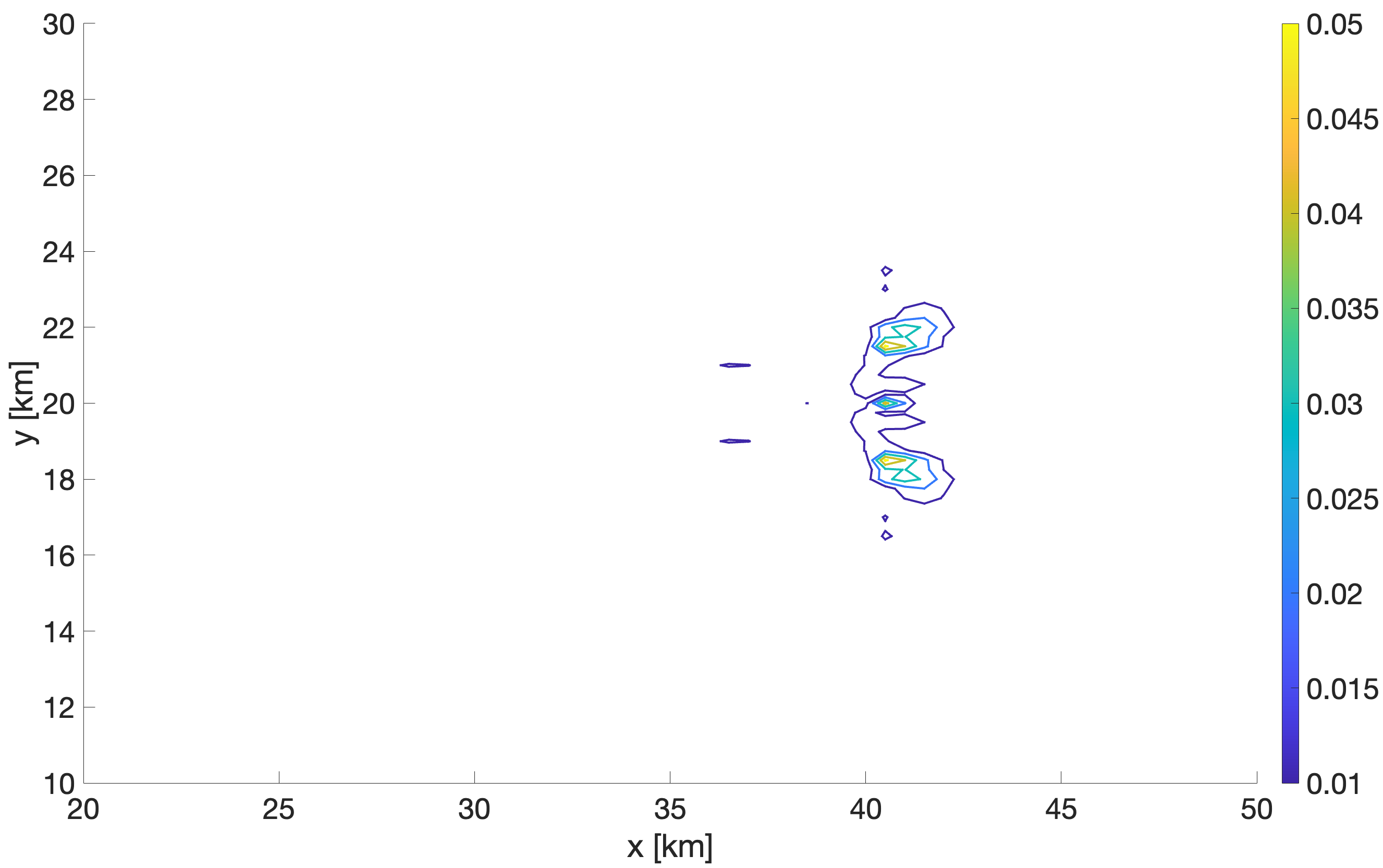}
	\end{subfigure}
	\\
	\caption{Vertical velocity in the 3D medium-steep bell-shaped hill test case at $T_{f} = \SI{10}{\hour}$, $x-y$ slice at $z = \SI{800}{\meter}$. Contours in the range $\SI[parse-numbers=false]{[-1.5, 1.3]}{\meter\per\second}$ with a $\SI{0.1}{\meter\per\second}$ interval. Top: comparison between the uniform mesh (black lines) and the non-conforming mesh (red lines). Negative contours are dashed. Bottom: absolute difference between the uniform mesh and the non-conforming mesh.}
	\label{fig:3D_contours_xy}
\end{figure}

\begin{figure}[h!]
	\centering
	\begin{subfigure}{0.75\textwidth}
		\centering
		\includegraphics[width = 0.925\textwidth]  {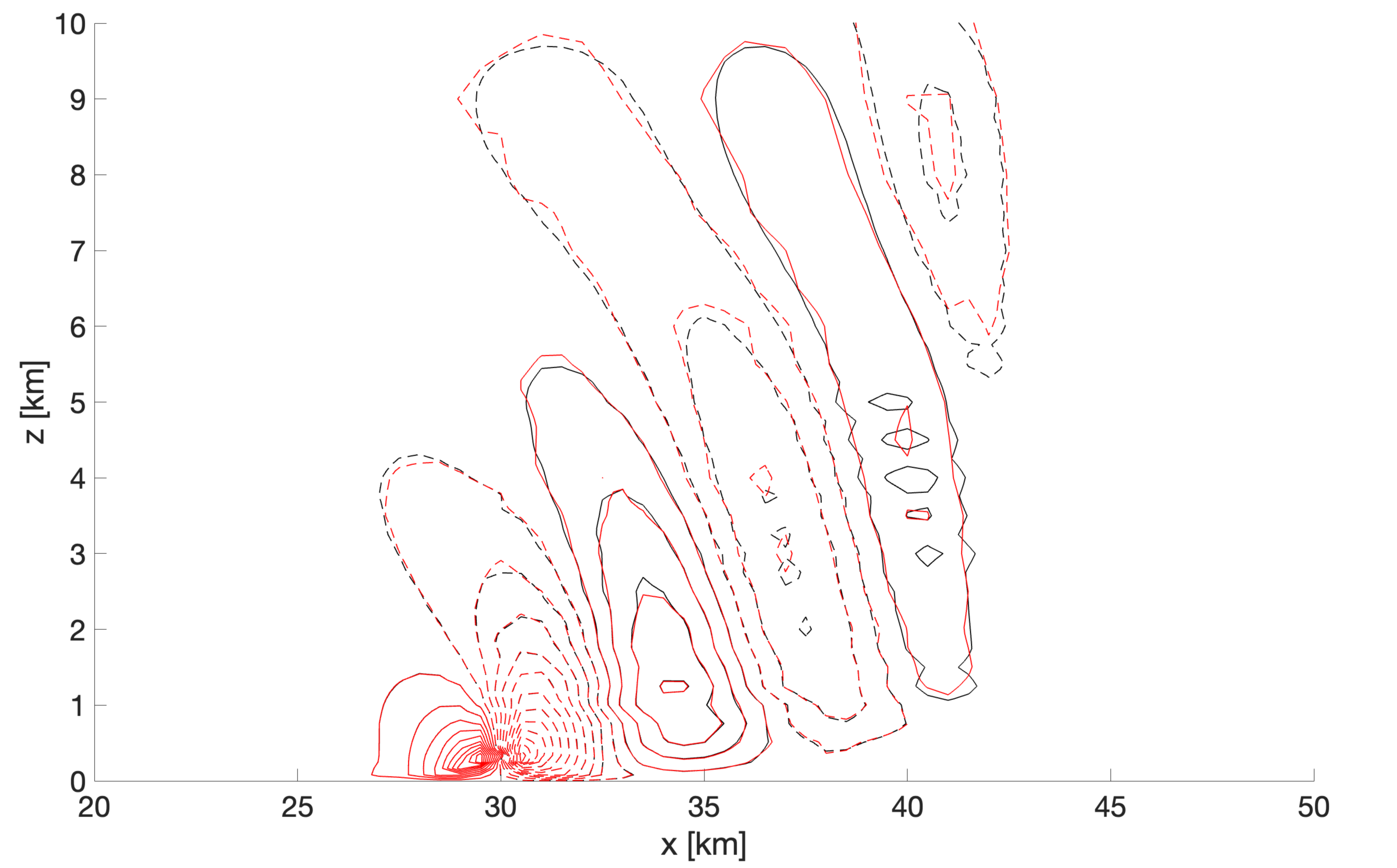}
	\end{subfigure}
	\\
	\begin{subfigure}{0.75\textwidth}
		\centering
		\includegraphics[width = 0.925\textwidth]{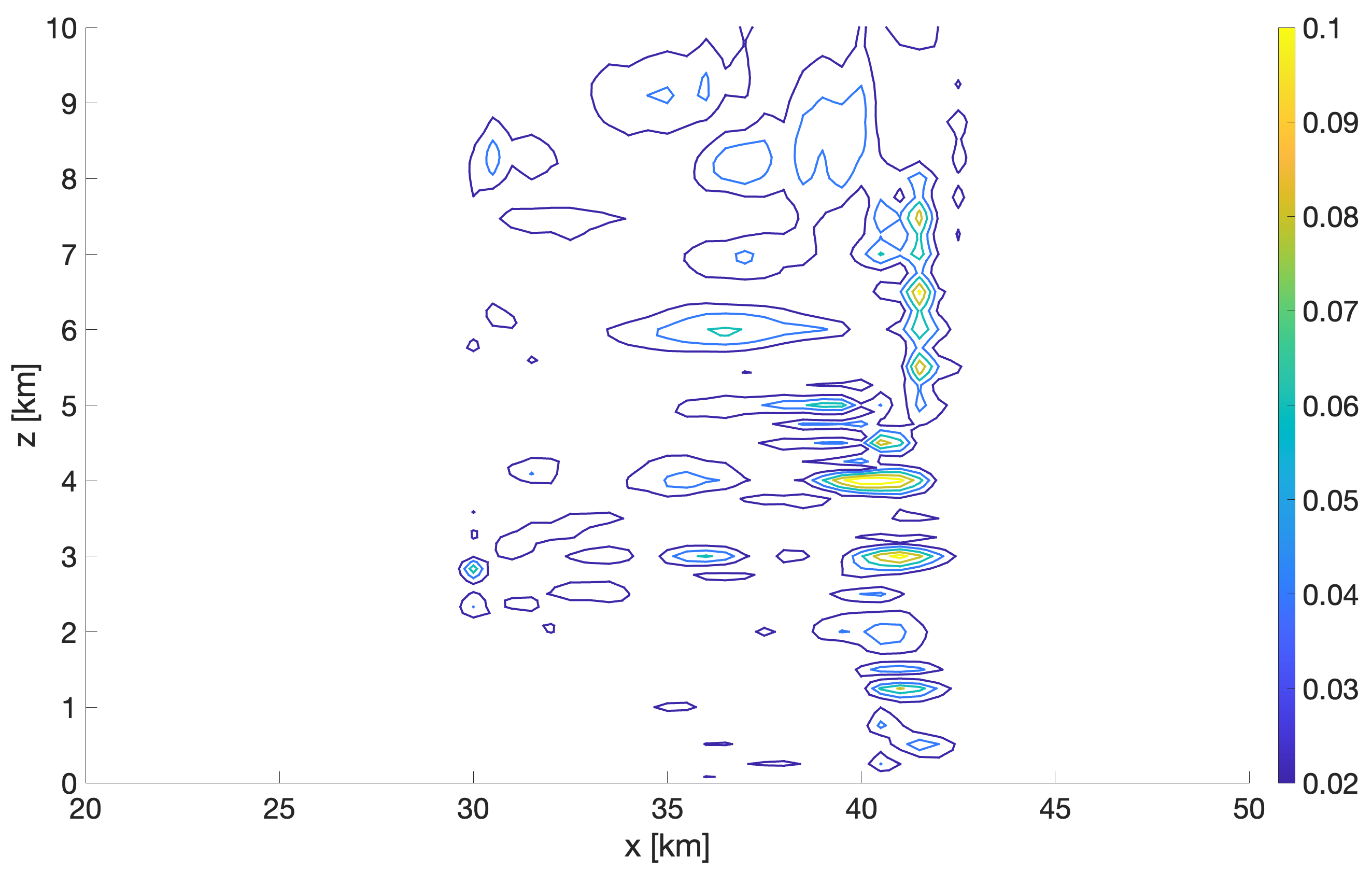}
	\end{subfigure}
	\\
	\caption{As in Figure \ref{fig:3D_contours_xy}, but showing an $x-z$ slice at $y = \SI{20}{\kilo\meter}$, and contours in the range $\SI[parse-numbers=false]{[-2.25, 2]}{\meter\per\second}$ with a $\SI{0.2}{\meter\per\second}$ interval.}
	\label{fig:3D_contours_xz}
\end{figure}

The contours plots of the vertical velocity on a $x-z$ slice placed at $y = \SI{20}{\kilo\meter}$ and on a $x-y$ slice at $z = \SI{800}{\meter}$ show the accuracy and the robustness of simulations employing non-conforming meshes, without significant differences compared to the simulation with uniform meshes (Figures \ref{fig:3D_contours_xy} and \ref{fig:3D_contours_xz}). Specifically, no spurious wave reflections arise at the internal boundaries that separate regions with different mesh resolutions. Hence, it is sufficient to employ a higher resolution only around the orography, whereas larger scales along all the directions can be resolved at a much coarser resolution.

\newpage

\subsection{Efficiency and scalability results}

Besides producing results that are of comparable accuracy with those obtained with a uniform mesh, the non-conforming mesh provides sizeable computational savings. With reference to the results in Figures \ref{fig:3D_contours_xy} and \ref{fig:3D_contours_xz}, the wall-clock time of the simulation with the non-conforming mesh is $\SI{2560}{\second}$, whereas the wall-clock time of the simulation with the uniform mesh is $\SI{36500}{\second}$. Hence, the use of the non-conforming mesh yields a computational time saving of around 93\% with respect to the uniform mesh (see also Table 6 in \cite{orlando:2024a}).

The size of this benchmark makes it a good candidate for a parallel scaling test. We consider a uniform mesh composed by $120 \times 80 \times 32 = 307200$ elements with polynomial degree $r = 4$, leading to around $38.5$ millions of unknowns for each scalar variable, and a non-conforming mesh composed by $N_{el} = 204816$ elements, yielding around $25.6$ millions of unknowns for each scalar variable. 

\begin{figure}[h!]
	\centering
	\begin{subfigure}{0.75\textwidth}
		\centering
		\includegraphics[width=0.9\textwidth]{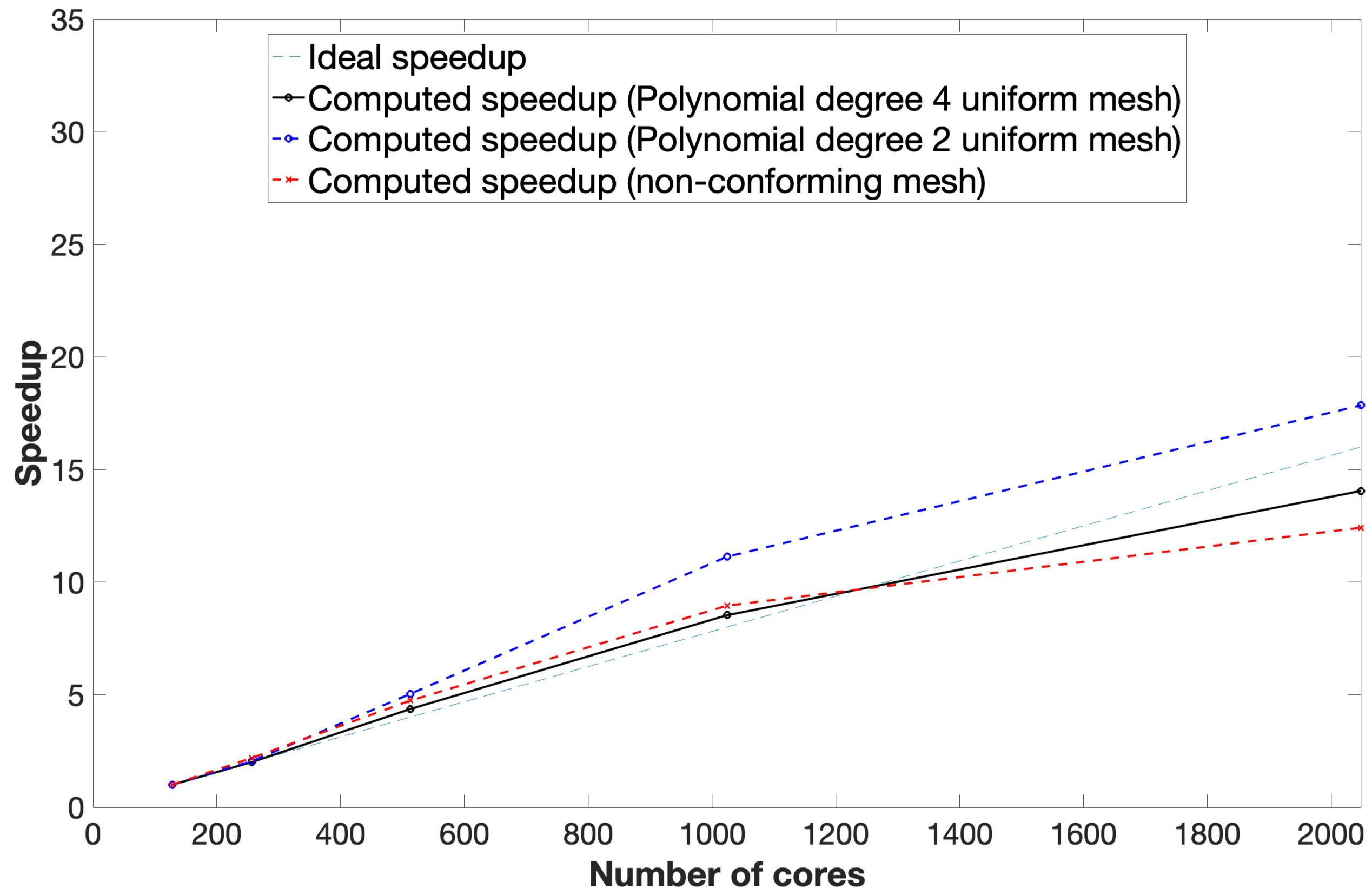}
	\end{subfigure}
	\\
	\begin{subfigure}{0.75\textwidth}
		\centering
		\includegraphics[width=0.9\textwidth]{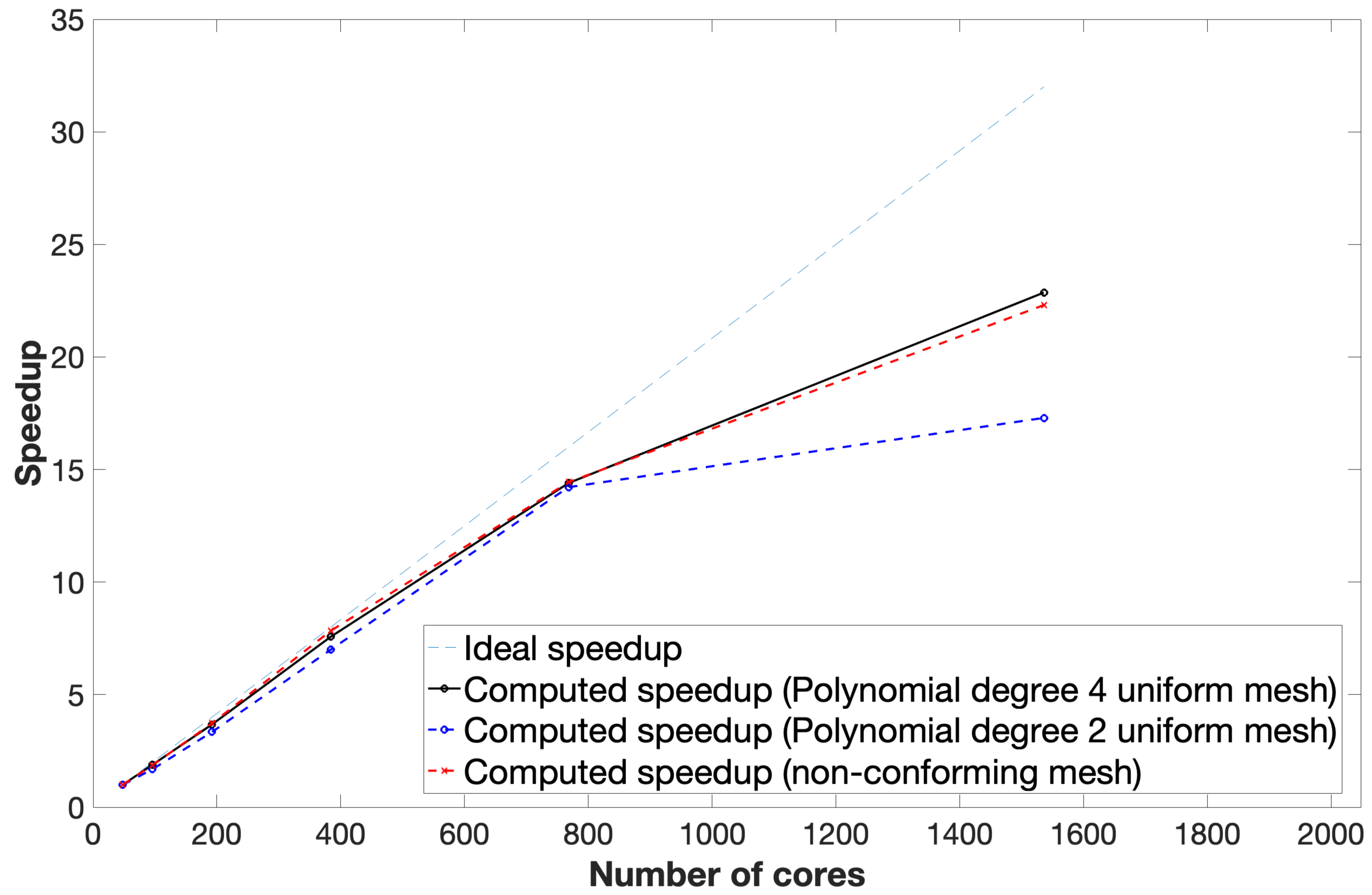}
	\end{subfigure}
	\caption{Strong scaling analysis. Computed speedup as a function of the number of cores used in the simulations with the uniform mesh using polynomial degree $4$ (solid black lines) and polynomial degree $2$ (dashed blue lines), and with the non-conforming mesh (dashed red lines). Top: exclusive use of computational nodes. Bottom: shared use of computational nodes.}
	\label{fig:strong_scaling}
\end{figure}

The strong scaling test evaluates the wallclock time of the simulations at fixed computational load (resolution) and using an increasing amount of computational resources. More specifically, we use 1, 2, 4, 8 and 16 full MeluXina CPU nodes with 128 cores each. In an ideal situation, the simulations speed up linearly with the amount of resources used. A good linear scaling is obtained up to 2048 cores, even with super-linear behaviour due to cache effects up to 1024 cores (Figure \ref{fig:strong_scaling}). In addition to the mentioned data locality of the discontinuous finite element approach, the favorable scaling is due to the matrix-free approach adopted in the solver, where no global sparse matrix is built and only the action of the linear operators is actually computed. These results represent a sensible improvement with respect to those previously presented in \cite{orlando:2023}, which evaluated scalability up to approximately half the cores used in this paper. The apparent better behaviour of the uniform mesh for a larger number of cores is due to the fact that more degrees of freedoms are involved and, therefore, the role of communication costs is less evident. 

In order to further emphasize this point, we perform an analogous scalability analysis with shared nodes, i.e. using computational nodes in which other jobs are simultaneously running. More specifically, we use 48, 96, 192, 384, 768 and 1536 MeluXina cores in shared nodes, and also perform runs at different polynomial degrees. One can easily notice that the use of shared nodes strongly degrades the parallel performance for more than about a thousand cores, and the effect is more marked at lower polynomial orders (Figure \ref{fig:strong_scaling}). Importantly, the speedup with the non-conforming mesh is comparable with that with the uniform mesh. Overall, the results confirm that the use of the Discontinuous Galerkin method, for which the stencil involves only the neighbours of each element, independently of the polynomial degree, provides an advantageous framework in terms of parallelization.

\begin{figure}[h!]
	\centering
	\begin{subfigure}{0.75\textwidth}
		\includegraphics[width=0.9\textwidth]{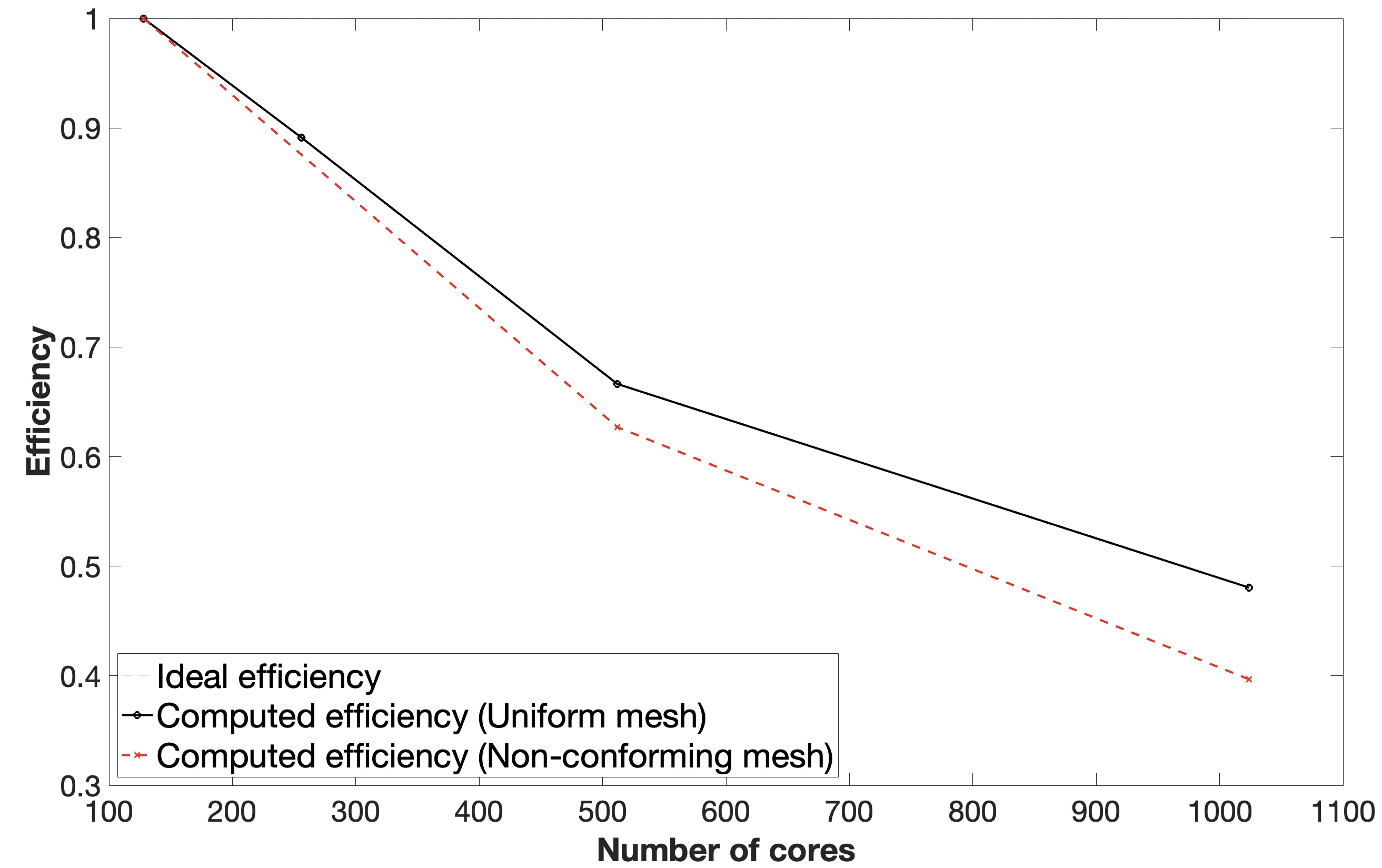}
	\end{subfigure}
	\begin{subfigure}{0.75\textwidth}
		\includegraphics[width=0.9\textwidth]{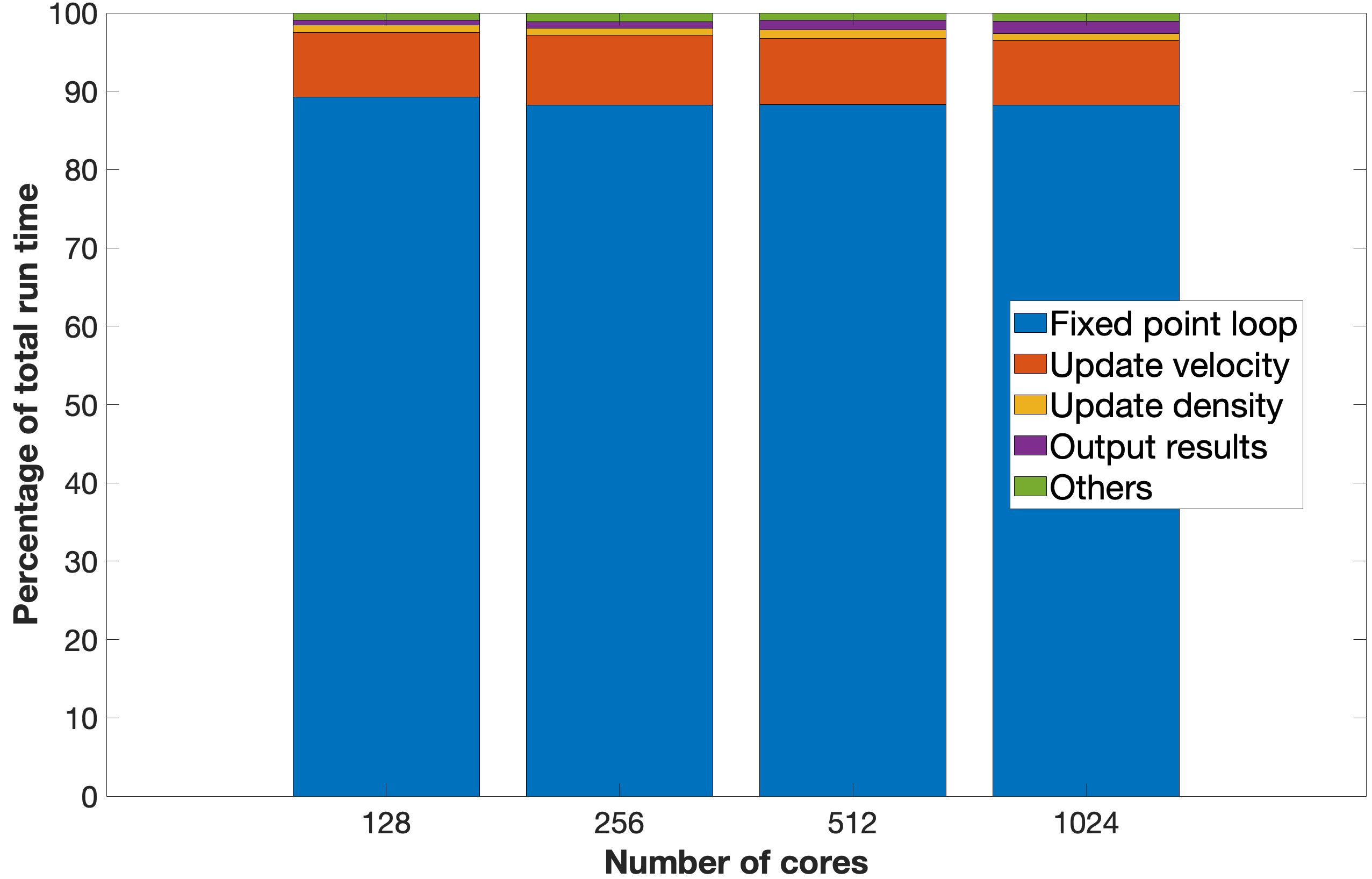}
	\end{subfigure}
	\caption{Weak scaling analysis. Top: parallel efficiency as a function of number of cores used in the simulations with the uniform mesh (solid black line) and the non-conforming mesh (dashed red line). Bottom: distribution of the computational time spent in various blocks of the algorithm as a function of number of cores.}
	\label{fig:weak_scaling}
\end{figure}

Moreover, a weak scaling analysis has been performed, using around $10^{5}$ unknowns per core for each scalar variable and increasing the problem size for an increasing amount of resources. For ideal scaling, wallclock time should not increase for increasing problem size. For the simulations both using the uniform mesh and using the non-conforming mesh, parallel performance is less than optimal in this case (Figure \ref{fig:weak_scaling}). However, the implementation actually outperforms the findings of previous \texttt{deal.II} studies \cite{kronbichler:2016}, where the efficiency of the Navier-Stokes solver implemented using the same library as in this paper drops to 20\%. In our simulations, a profiling study reveals that most of the time is spent in the fixed point loop to update the pressure variable, for which a non-symmetric linear system arises and a GMRES solver is therefore employed \cite{orlando:2022, orlando:2023, orlando:2024a}. However, in terms of percentage time with respect to the total run time, the time spent in the linear solver is similar with increasing core counts (Figure \ref{fig:weak_scaling}). It is therefore expected that an improved solver strategy based on, e.g., advanced preconditioning techniques, will improve efficiency independently of the computational resources employed.

\clearpage

\section{Conclusions}
\label{sec:conclu}

This paper has reported results of performance tests with a new high-order Discontinuous Galerkin model for atmospheric dynamics simulations using non-conforming mesh refinement. On a three-dimensional Cartesian benchmark of dry compressible fluid flow over orography with a stably stratified background atmosphere, the simulations using non-conforming meshes provide results that are equally accurate and more than 90\% more efficient than simulations using uniform meshes that are standard in operational atmospheric modelling. In addition, the data locality features of the matrix-free, discontinuous finite element-based approach ensure good CPU scalability as measured in parallel runs with the state-of-the-art MeluXina EuroHPC facility. 

These results open up a number of future avenues for investigations. First, enhanced realistic simulations will come from the  inclusion of more complex physical phenomena, in particular moist air, and the use of more general equations of state for real gases, for which the numerical method proposed in \cite{orlando:2022} has been already shown to be effective.
Next, the development of a three-dimensional dynamical core in spherical geometry including rotation will enable the testing on more realistic atmospheric flows, for which more sizeable computational resources will be required. This more general and computationally heavier context will make it easier to fine-tune the performance and improve the findings in this paper, especially regarding weak scaling, and to gauge the viability of the proposed model towards full-fledged numerical weather prediction capability.

\section*{Acknowledgements}

The simulations have been run thanks to the computational resources made available through the EuroHPC JU Benchmark And Development project EHPC-BEN-2024B03-045. This work has been partly supported by the ESCAPE-2 project, European Union’s Horizon 2020 Research and Innovation Programme (Grant Agreement No. 800897).

\printbibliography

\end{document}